\documentclass[11pt,a4paper,reqno]{amsart}  
\usepackage{amssymb} \usepackage{amscd} \usepackage{times}
\newcommand{\bea} {\begin{eqnarray}}
\newcommand{\eea} {\end{eqnarray}}
\newcommand{\Bea} {\begin{eqnarray*}}
\newcommand{\Eea} {\end{eqnarray*}}

\def\zbb{\mathbb{Z}}  
  
  \def\phi{\varphi}
 \def\p1{{\mathbb{P}^1_\zbb}}

\newtheorem{Theorem}{\quad Theorem}[section]

\usepackage{amsmath,  amsfonts}

\usepackage[T1]{fontenc}

\newcommand{\be} {\begin{equation}}
\newcommand{\ee} {\end{equation}}

\begin{document}
\title{About Brezis-Merle Problem with Lipschitz condition.}
\author{Samy Skander Bahoura} 
\address{Departement de Mathematiques, Universite Pierre et Marie Curie, 4 place Jussieu, 75005, Paris, France.}
\email{samybahoura@gmail.com, samybahoura@yahoo.fr}
\maketitle
\begin{abstract}
We give blow-up analysis for a Brezis and Merle's problem with Dirichlet condition. As an application we have a proof of a compactness result under Lipschitz condition on the prescribed scalar curvature and a weaker assumption on the regularity of the domain (smooth domain or $ C^{2,\alpha} $ domain, $ 1\geq \alpha >0 $). \end{abstract}


\section{Introduction and Main Results} 

We set $ \Delta = -(\partial_{11} +\partial_{22}) $  on open set $ \Omega $ of $ {\mathbb R}^2 $ with a smooth (or $ C^{2,\alpha}, \alpha >0 $) boundary.

We consider the following equation:

\begin{displaymath}  (P)  \left \{ \begin {split} 
      \Delta u  & = V e^{u}     \,\, &&\text{in} \!\!&&\Omega \subset {\mathbb R}^2, \\
                  u  & = 0  \,\,             && \text{in} \!\!&&\partial \Omega.               
\end {split}\right.
\end{displaymath}

Here, we assume that:

$$ 0 \leq V \leq b < + \infty, \,\, e^u \in L^1({\Omega})\,\, {\rm and} \,\,  u \in W_0^{1,1}(\Omega). $$

We can see in [7] a nice formulation of this problem $ (P) $  in the sense of the distributions. This Problem arises from geometrical and physical problems see for example [1, 2, 18, 19]. The above equation was studied by many authors, with or without  the boundary condition, also for Riemannian surfaces,  see [1-19],  where one can find some existence and compactness results. In [6] we have the following important Theorem,

{\bf Theorem A}{\it (Brezis-Merle [6])}.{\it For $ (u_i)_i $ and $ (V_i)_i $ two sequences of functions relative to $ (P) $ with,

$$ 0 < a \leq V_i \leq b < + \infty $$

then, for all compact subset $ K $ of $\Omega $ it holds,

$$ \sup_K u_i \leq c, $$

with $ c $ depending on $ a, b,  K $ and $ \Omega $.}

One can find in [6] an interior estimate if we assume $ a=0 $, but we need an assumption on the integral of $ e^{u_i} $, namely, we have:

{\bf Theorem B}{\it (Brezis-Merle [6])}.{\it For $ (u_i)_i $ and $ (V_i)_i $ two sequences of functions relative to the problem $ (P) $ with,

$$ 0 \leq V_i \leq b < + \infty \,\, {\rm and} \,\, \int_{\Omega} e^{u_i} dy  \leq C, $$

then, for all compact subset $ K $ of $ \Omega $ it holds;

$$ \sup_K u_i \leq c, $$

with $ c $ depending on $ b, C, K $ and $ \Omega $.}

We look to the uniform boundedness on all $ \bar \Omega  $ of sequences of solutions of the Problem $ (P) $. Remark that, when $ a=0 $ the boundedness of $ \int_{\Omega} e^{u_i} $ is a necessary condition in the problem $ (P) $ as showed in $ [6] $ by the following counterexample.

{\bf Theorem C}{\it (Brezis-Merle [6])}.{\it There are two sequences $ (u_i)_i $ and $ (V_i)_i $ of the problem $ (P) $ with, 

$$ 0 \leq V_i \leq b < + \infty \,\, {\rm and} \,\, \int_{\Omega} e^{u_i} dy  \leq C, $$

such that,

$$ \sup_{\Omega}  u_i \to + \infty. $$}

To obtain the two first previous results (Theorems A and B) Brezis and Merle used  an inequality (Theorem 1 of [6]) obtained by an approximation argument with the Fatou's lemma and they applied the maximum principle in $ W_0^{1,1}(\Omega) $ which arises from Kato's inequality. Also this weak form of the maximum principle is used to prove the local uniform boundedness result by comparing  a certain function and the Newtonian potential. We refer to [5] for a topic about the weak form of the maximum principle.

\smallskip

{\bf Remarks:} 1) Theorem 1 of [6], can be obtained by the usual maximum principle and Agmon regularity theorem which require $ C^2 $ regularity on the domain.

\smallskip

2) The duality Theorem which we use require $ C^2 $ regularity on the domain, see Gilbarg-Trudinger books.

\smallskip

Note that for the problem $ (P) $, by using the Pohozaev identity, we can prove that $ \int_{\Omega} e^{u_i} $ is uniformly bounded when $ 0 < a \leq V_i \leq b < +\infty $ and $  ||\nabla V_i||_{L^{\infty}} \leq A $ and $ \Omega $ starshaped, when $ a=0 $ and $ \nabla \log V_i $ is uniformly bounded, we can bound uniformly $ \int_{\Omega} V_i e^{u_i} $. In [16] Ma-Wei have proved that those results stay true for all open sets not necessarily starshaped.

In [9] Chen-Li have proved that if $ a=0 $ and $ \nabla \log V_i $ is uniformly bounded, then  the functions are uniformly bounded near the boundary.

In [9] Chen-Li have proved that if $ a=0 $ and $ \int_{\Omega} e^{u_i} $ is uniformly bounded and $ \nabla \log V_i $ is uniformly bounded, then we have the compactness result directly. Ma-Wei in [16], extend this result in the case where $ a >0 $.

If we assume $ V $ more regular, we can have another type of estimates called $ \sup + \inf $ type inequalities. It was proved by Shafrir see [17] that, if $ (u_i)_i, (V_i)_i $ are two sequences of functions solutions of the previous equation without assumption on the boundary and, $ 0 < a \leq V_i \leq b < + \infty $, then we have the following interior estimate:

$$ C\left (\dfrac{a}{b} \right ) \sup_K u_i + \inf_{\Omega} u_i \leq c=c(a, b, K, \Omega). $$

One can see in [10] an explicit value of $ C\left (\dfrac{a}{b}\right ) =\sqrt {\dfrac{a}{b}} $. In his proof Shafrir has used a blow-up function, the Stokes formula and an isoperimetric inequality see [2]. For Chen-Lin, they have used the blow-up analysis combined with some geometric type inequality for the integral curvature.

\smallskip

Now, if we suppose $ (V_i)_i $ uniformly Lipschitzian with $ A $ its
Lipschitz constant then $ C(a/b)=1 $ and $ c=c(a, b, A, K, \Omega)
$ see Brezis-Li-Shafrir [4]. This result was extended for
H\"olderian sequences $ (V_i)_i $ by Chen-Lin, see  [10]. Also, one
can see in [14] an extension of the Brezis-Li-Shafrir result to compact
Riemannian surfaces without boundary. One can see in [15] explicit form,
($ 8 \pi m, m\in {\mathbb N}^* $ exactly), for the numbers in front of
the Dirac masses when the solutions blow-up. Here, the notion of isolated blow-up point is used.

 In [8] we have some a priori estimates on the 2 and 3-spheres $ {\mathbb S}_2 $, $ {\mathbb S}_3 $.

Here we give the behavior of the blow-up points on the boundary and a proof of Brezis-Merle Problem with Lipschitz condition.

The Brezis-Merle Problem (see [6]) is:

{\bf Problem}. Suppose that $ V_i \to  V $ in $ C^0( \bar \Omega ) $ with $ 0 \leq V_i $. Also, we consider a sequence of solutions $ (u_i) $ of $ (P) $ relative to $ (V_i) $ such that,

$$ \int_{\Omega} e^{u_i} dx \leq C,  $$

is it possible to have:

$$ ||u_i||_{L^{\infty}}\leq C ? $$

Here, we give a caracterization of the behavior of the blow-up points on the boundary and also, in particular we extend Chen-Li theorems, indeed, the result of Chen-Li holds for analytic domains and our result holds for smooth of $ C^{2,\alpha} $ domains. For the behavior of the blow-up points on the boundary, the following condition is enough,

$$ 0 \leq  V_i \leq b. $$

The condition $ V_i \to  V $ in $ C^0(\bar \Omega) $ is not necessary, but for the proof of the compactness for the Brezis-Merle problem we assume that:

$$ ||\nabla V_i||_{L^{\infty}}\leq  A. $$

Our main results are:

\begin{Theorem}  Assume that $ \max_{\Omega} u_i \to +\infty $, where $ (u_i) $ are solutions of the problem $ (P) $ with:
 
 $$ 0 \leq V_i \leq b\,\,\, {\rm and } \,\,\, \int_{\Omega}  e^{u_i} dx \leq C, \,\,\, \forall \,\, i, $$
 
 then, after passing to a subsequence, there is a function $ u $,  there is a number $ N \in {\mathbb N} $ and there are $ N  $ points $ x_1, \ldots, x_N \in  \partial \Omega $, such that, 

$$ \partial_{\nu} u_i  \to \partial_{\nu} u +\sum_{j=1}^N \alpha_j \delta_{x_j}, \,\,\, \alpha_j \geq 4\pi, \,\, {\rm in \, the \, sense \, of \, measures \, on \,\, \partial \Omega.} $$

$$ u_i \to u \,\,\, {\rm in }\,\,\, C^1_{loc}(\bar \Omega-\{x_1,\ldots, x_N \}). $$
 
\end{Theorem} 

\begin{Theorem}Assume that $ (u_i) $ are solutions of $ (P) $ relative to $ (V_i) $ with the following conditions:

$$ 0 \leq V_i \leq b, \,\, ||\nabla V_i||_{L^{\infty}} \leq A \,\,\, {\rm and } \,\,\, \int_{\Omega} e^{u_i} \leq C, $$

we have,

$$  || u_i||_{L^{\infty}} \leq c(b, A, C, \Omega). $$

\end{Theorem}
 
In the previous theorem we have a proof of the global a priori estimate which concern the problem $ (P) $. The proof of Chen-Li and Ma-Wei [9,16], use the moving-plane method for the case $ \nabla \log V_i $ uniformly bounded near the boundary (and $ C^{2,\alpha} $ domain, $ 1\geq \alpha >0 $) and for analytic domain for the case $ \nabla V_i $ uniformly bounded.

\smallskip

To prove Theorem 1.2, we argue by contradiction and use Theorem 1.1.

\section{Proof of the theorems} 

\underbar {\it Proof of theorem 1.1:} 

We have:

$$ u_i \in W_0^{1,1}(\Omega). $$

Since $ e^{u_i} \in L^1(\Omega) $ by the corollary 1 of Brezis-Merle's paper (see [6]) we have $ e^{u_i} \in L^k(\Omega) $ for all $ k  >2 $ and the elliptic estimates of Agmon and the Sobolev embedding (see [1]) imply that:

$$ u_i \in W^{2, k}(\Omega)\cap C^{1, \epsilon}(\bar \Omega). $$ 

We denote by $ \partial_{\nu} u_i $ the inner normal derivative. By the maximum principle we have, $ \partial_{\nu} u_i  \geq 0 $.

By the Stokes formula we have, 

$$ \int_{\partial \Omega} \partial_{\nu} u_i d\sigma \leq C, $$

We use the weak convergence in the space of Radon measures to have the existence of a nonnegative Radon measure $ \mu $ such that,

$$ \int_{\partial \Omega} \partial_{\nu} u_i \phi  d\sigma \to \mu(\phi), \,\,\, \forall \,\,\, \phi \in C^0(\partial \Omega). $$

We take an $ x_0 \in \partial \Omega $ such that, $ \mu({x_0}) < 4\pi $. For $ \epsilon >0 $ small enough set $ I_{\epsilon}= B(x_0, \epsilon)\cap \partial \Omega $. We choose a function $ \eta_{\epsilon} $ such that,

$$ \begin{cases}
    
\eta_{\epsilon} \equiv 1,\,\,\,  {\rm on } \,\,\,  I_{\epsilon}, \,\,\, 0 < \epsilon < \delta/2,\\

\eta_{\epsilon} \equiv 0,\,\,\, {\rm outside} \,\,\, I_{2\epsilon }, \\

0 \leq \eta_{\epsilon} \leq 1, \\

||\nabla \eta_{\epsilon}||_{L^{\infty}(I_{2\epsilon})} \leq \dfrac{C_0(\Omega, x_0)}{\epsilon}.

\end{cases} $$

We take a $\tilde \eta_{\epsilon} $ such that,

\begin{displaymath}  \left \{ \begin {split} 
      \Delta \tilde \eta_{\epsilon}  & = 0              \,\, &&\text{in} \!\!&&\Omega \subset {\mathbb R}^2, \\
                  \tilde\eta_{\epsilon} & =  \eta_{\epsilon}   \,\,             && \text{in} \!\!&&\partial \Omega.               
\end {split}\right.
\end{displaymath}

{\bf Remark:} We use the following steps in the construction of $ \eta_{\epsilon} $:

We take a cutoff function $ \eta_{0} $ in $ B(0, 2) $ or $ B(x_0, 2) $:

1- We set $ \eta_{\epsilon}(x)= \eta_0(|x-x_0|/\epsilon) $ in the case of the unit disk it is sufficient.

2- Or, in the general case: we use a chart $ (f, \tilde \Omega =f(B_r(0))) $, for $ r>0 $ small enough and $ f(0)=x_0 $ and we take $ \mu_{\epsilon}(x)= \eta_0 ( f( |x|/ \epsilon)) $ to have  connected  sets $ I_{\epsilon} $ and we take $ \eta_{\epsilon}(y)= \mu_{\epsilon}(f^{-1}(y))$. Because $ f, f^{-1} $ are Lipschitz, $ |f(x)-x_0| \leq k_ 2|x|\leq 1 $ for $ |x| \leq 1/k_2 $ and $ |f(x)-x_0| \geq k_ 1|x|\geq 2 $ for $ |x| \geq 2/k_1>1/k_2 $, the support  of $ \eta $ is in $ I_{(2/k_1)\epsilon} $.

$$ \begin{cases}
    
\eta_{\epsilon} \equiv 1,\,\,\,  {\rm on } \,\,\,  f(I_{(1/k_2)\epsilon}), \,\,\, 0 < \epsilon < \delta/2,\\

\eta_{\epsilon} \equiv 0,\,\,\, {\rm outside} \,\,\, f(I_{(2/k_1)\epsilon }), \\

0 \leq \eta_{\epsilon} \leq 1, \\

||\nabla \eta_{\epsilon}||_{L^{\infty}(I_{(2/k_1)\epsilon})} \leq \dfrac{C_0(\Omega, x_0)}{\epsilon}.

\end{cases} $$

3- Also, we can take: $ \mu_{\epsilon}(x)= \eta_0(|x|/\epsilon) $ and $ \eta_{\epsilon}(y)= \mu_{\epsilon}(f^{-1}(y)) $, we extend it by $ 0 $ outside $ f(B_1(0)) $.  We have $ f(B_1(0)) = D_1(x_0) $, $ f (B_{\epsilon}(0))= D_{\epsilon}(x_0) $ and $ f(B_{\epsilon}^+)= D_{\epsilon}^+(x_0) $ with $ f $ and $ f^{-1} $ smooth diffeomorphism.

$$ \begin{cases}
    
\eta_{\epsilon} \equiv 1,\,\,\,  {\rm on \, the \, connected \, set } \,\,\,  J_{\epsilon} =f(I_{\epsilon}), \,\,\, 0 < \epsilon < \delta/2,\\

\eta_{\epsilon} \equiv 0,\,\,\, {\rm outside} \,\,\, J'_{\epsilon}=f(I_{2\epsilon }), \\

0 \leq \eta_{\epsilon} \leq 1, \\

||\nabla \eta_{\epsilon}||_{L^{\infty}(J'_{\epsilon})} \leq \dfrac{C_0(\Omega, x_0)}{\epsilon}.

\end{cases} $$

And, $ H_1(J'_{\epsilon}) \leq C_1 H_1(I_{2\epsilon}) = C_1 4\epsilon $, because $ f $ is Lipschitz. Here $ H_1 $ is the Hausdorff measure.

We solve the Dirichlet Problem:

\begin{displaymath}  \left \{ \begin {split} 
      \Delta \bar \eta_{\epsilon}  & = \Delta \eta_{\epsilon}              \,\, &&\text{in} \!\!&&\Omega \subset {\mathbb R}^2, \\
                  \bar \eta_{\epsilon} & = 0   \,\,             && \text{in} \!\!&&\partial \Omega.               
\end {split}\right.
\end{displaymath}

and finaly we set $ \tilde \eta_{\epsilon} =-\bar \eta_{\epsilon} + \eta_{\epsilon} $. Also, by the maximum principle and the elliptic estimates we have :

$$ ||\nabla \tilde \eta_{\epsilon}||_{L^{\infty}} \leq C(|| \eta_{\epsilon}||_{L^{\infty}} +||\nabla \eta_{\epsilon}||_{L^{\infty}} + ||\Delta \eta_{\epsilon}||_{L^{\infty}}) \leq \dfrac{C_1}{\epsilon^2}, $$

with $ C_1 $ depends on $ \Omega $.

We use the following estimate, see [3, 7, 19],

$$ ||\nabla u_i||_{L^q} \leq C_q, \,\,\forall \,\, i\,\, {\rm and  }  \,\, 1< q < 2. $$

We deduce from the last estimate that, $ (u_i) $ converge weakly in $ W_0^{1, q}(\Omega) $, almost everywhere to a function $ u \geq 0 $ and $ \int_{\Omega} e^u < + \infty $ (by Fatou's lemma). Also, $ V_i $ weakly converge to a nonnegative function $ V $ in $ L^{\infty} $. The function $ u $ is in $ W_0^{1, q}(\Omega) $ solution of :

\begin{displaymath} \left \{ \begin {split} 
      \Delta u  & = V e^{u} \in L^1(\Omega)    \,\, &&\text{in} \!\!&&\Omega \subset {\mathbb R}^2, \\
                  u  & = 0  \,\,                                     && \text{in} \!\!&&\partial \Omega.               
\end {split}\right.
\end{displaymath}
 
According to the corollary 1 of Brezis-Merle result, see [6],   we have $ e^{k u }\in L^1(\Omega), k >1 $. By the elliptic estimates, we have $ u \in C^1(\bar \Omega) $.

For two vectors $ v, w $ of $ {\mathbb R}^2 $ we denote by $ v \cdot w $ the inner product of $ v $ and $ w $.

We can write,

\be \Delta ((u_i-u) \tilde \eta_{\epsilon})= (V_i e^{u_i} -Ve^u)\tilde \eta_{\epsilon} -2\nabla (u_i- u)\cdot \nabla \tilde \eta_{\epsilon} . \label{(1)}\ee

We use the interior esimate of Brezis-Merle, see [6],

\underbar {\it Step 1:} Estimate of the integral of the first term of the right hand side of $ (\ref{(1)}) $.

We use the Green formula between $ \tilde \eta_{\epsilon} $ and $ u $, we obtain,

\be  \int_{\Omega} Ve^u \tilde \eta_{\epsilon} dx =\int_{\partial \Omega} \partial_{\nu} u \eta_{\epsilon} \leq C'\epsilon ||\partial_{\nu}u||_{L^{\infty}}= C \epsilon \label{(2)}\ee

We have,

\begin{displaymath} \left \{ \begin {split} 
      \Delta u_i  & = V_i e^{u_i}                      \,\, &&\text{in} \!\!&&\Omega \subset {\mathbb R}^2, \\
                  u_i  & = 0  \,\,                                     && \text{in} \!\!&&\partial \Omega.               
\end {split}\right.
\end{displaymath}

We use the Green formula between $ u_i $ and $ \tilde \eta_{\epsilon} $ to have:

\be \int_{\Omega} V_i e^{u_i} \tilde \eta_{\epsilon} dx = \int_{\partial \Omega} \partial_{\nu} u_i \eta_{\epsilon} d\sigma \to \mu(\eta_{\epsilon}) \leq \mu(J'_{\epsilon}) \leq 4  \pi - \epsilon_0, \,\,\, \epsilon_0 >0 \label{(3)}\ee

From $ (\ref{(2)}) $ and $ (\ref{(3)}) $ we have for all $ \epsilon >0 $ there is $ i_0 =i_0(\epsilon) $ such that, for $ i \geq i_0 $,

\be \int_{\Omega} |(V_ie^{u_i}-Ve^u) \tilde \eta_{\epsilon}| dx \leq 4 \pi -\epsilon_0 + C \epsilon \label{(4)}\ee

\underbar {\it Step 2:} Estimate of integral of the second term of the right hand side of $ (\ref{(1)}) $.

Let $ \Sigma_{\epsilon} = \{ x \in \Omega, d(x, \partial \Omega) = \epsilon^3 \} $ and $ \Omega_{\epsilon^3} = \{ x \in \Omega, d(x, \partial \Omega) \geq \epsilon^3 \} $, $ \epsilon > 0 $. Then, for $ \epsilon $ small enough, $ \Sigma_{\epsilon} $ is a manifold.

The measure of $ \Omega-\Omega_{\epsilon^3} $ is $ k_2\epsilon^3 \leq meas(\Omega-\Omega_{\epsilon^3}) = \mu_L (\Omega-\Omega_{\epsilon^3}) \leq k_1 \epsilon^3 $. Here $ \mu_L $ is the Lebesgue measure.

{\bf Remark}: for the unit ball $ {\bar B(0,1)} $, our new manifold is $ {\bar B(0, 1-\epsilon^3)} $.

( Proof of this fact; let's consider $ d(x, \partial \Omega) = d(x, z_0), z_0 \in \partial \Omega $, this imply that $ (d(x,z_0))^2 \leq (d(x, z))^2 $ for all $ z \in \partial \Omega $ which it is equivalent to $ (z-z_0) \cdot (2x-z-z_0) \leq 0 $ for all $ z \in \partial \Omega $, let's consider a chart around $ z_0 $ and $ \gamma (t) $ a curve in $ \partial \Omega $, we have;

$ (\gamma (t)-\gamma(t_0) \cdot (2x-\gamma(t)-\gamma(t_0)) \leq 0 $ and it is clear that, $ \gamma'(t_0) \cdot (x-\gamma(t_0)) = 0 $, which imply that $ x= z_0-s \nu_0 $ where $ \nu_0 $ is the outward normal of $ \partial \Omega $ at $ z_0 $))

With this fact, we can say that $ S= \{ x, d(x, \partial \Omega) \leq \epsilon \}= \{ x= z_0- s \nu_{z_0}, z_0 \in \partial \Omega, \,\, -\epsilon \leq s \leq \epsilon \} $. It  is sufficient to work on  $ \partial \Omega $. Let's consider a charts $ (z, D=B(z, 4 \epsilon_z), \gamma_z) $ with $ z \in \partial \Omega $ such that $ \cup_z B(z, \epsilon_z) $ is  cover of $ \partial \Omega $ .  One can extract a finite cover $ (B(z_k, \epsilon_k)), k =1, ..., m $, by the area formula the measure of $ S \cap B(z_k, \epsilon_k) $ is less than a $ k\epsilon $ (a $ \epsilon $-rectangle). For the reverse inequality, it is sufficient to consider one chart around one point on the boundary). 

We write,

\be \int_{\Omega} |\nabla ( u_i -u) \cdot \nabla \tilde \eta_{\epsilon}| dx =
\int_{\Omega_{\epsilon^3}} |\nabla (u_i -u) \cdot \nabla \tilde \eta_{\epsilon}| dx + \int_{\Omega - \Omega_{\epsilon^3}} |\nabla (u_i-u) \cdot \nabla \tilde \eta_{\epsilon}| dx.  \label{(5)}\ee

\underbar {\it Step 2.1:} Estimate of $ \int_{\Omega - \Omega_{\epsilon^3}} |\nabla (u_i-u) \cdot \nabla \tilde \eta_{\epsilon}| dx $.

First, we know from the elliptic estimates that  $ ||\nabla \tilde \eta_{\epsilon}||_{L^{\infty}} \leq C_1 /\epsilon^2 $, $ C_1 $ depends on $ \Omega $

We know that $ (|\nabla u_i|)_i $ is bounded in $ L^q, 1< q < 2 $, we can extract  from this sequence a subsequence which converge weakly to $ h \in L^q $. But, we know that we have locally the uniform convergence to $ |\nabla u| $ (by Brezis-Merle's theorem), then, $ h= |\nabla u| $ a.e. Let $ q' $ be the conjugate of $ q $.

We have, $  \forall f \in L^{q'}(\Omega)$

$$ \int_{\Omega} |\nabla u_i| f dx \to \int_{\Omega} |\nabla u| f dx $$

If we take $ f= 1_{\Omega-\Omega_{\epsilon^3}} $, we have:

$$ {\rm for } \,\, \epsilon >0 \,\, \exists \,\, i_1 = i_1(\epsilon) \in {\mathbb N}, \,\,\, i \geq  i_1,  \,\, \int_{\Omega-\Omega_{\epsilon^3}} |\nabla u_i| \leq \int_{\Omega-\Omega_{\epsilon^3}} |\nabla u| + \epsilon^3. $$

Then, for $ i \geq i_1(\epsilon) $,

$$ \int_{\Omega-\Omega_{\epsilon^3}} |\nabla u_i| \leq meas(\Omega-\Omega_{\epsilon^3}) ||\nabla u||_{L^{\infty}} + \epsilon^3 = \epsilon^3(k_1 ||\nabla u||_{L^{\infty}} + 1). $$

Thus, we obtain,

\be \int_{\Omega - \Omega_{\epsilon^3}} |\nabla (u_i-u) \cdot \nabla \tilde \eta_{\epsilon}| dx \leq  \epsilon C_1(2 k_1 ||\nabla u||_{L^{\infty}} + 1) \label{(6)}\ee

The constant $ C_1 $ does  not depend on $ \epsilon $ but on $ \Omega $.

\underbar {\it Step 2.2:} Estimate of $ \int_{\Omega_{\epsilon^3}} |\nabla (u_i-u) \cdot \nabla \tilde \eta_{\epsilon}| dx $.

We know that, $ \Omega_{\epsilon} \subset \subset \Omega $, and ( because of Brezis-Merle's interior estimates) $ u_i \to u $ in $ C^1(\Omega_{\epsilon^3}) $. We have,

$$ ||\nabla (u_i-u)||_{L^{\infty}(\Omega_{\epsilon^3})} \leq \epsilon^3,\, {\rm for } \,\, i \geq i_3 = i_3(\epsilon). $$

We write,
 
$$ \int_{\Omega_{\epsilon3}} |\nabla (u_i-u) \cdot \nabla \tilde \eta_{\epsilon} | dx \leq ||\nabla (u_i-u)||_{L^{\infty}(\Omega_{\epsilon^3})} ||\nabla \tilde \eta_{ \epsilon}||_{L^{\infty}} \leq C_1 \epsilon \,\, {\rm for } \,\, i \geq i_3, $$

For $ \epsilon >0 $, we have for $ i \in {\mathbb N} $, $ i \geq \max \{i_1, i_2, i_3 \} $,

\be \int_{\Omega} |\nabla (u_i-u) \cdot \nabla \tilde \eta_{\epsilon}| dx \leq \epsilon C_1(2 k_1 ||\nabla u||_{L^{\infty}} + 2) \label{(7)}\ee

From $ (\ref{(4)}) $ and $ (\ref{(7)}) $, we have, for $ \epsilon >0 $, there is $ i_3= i_3(\epsilon) \in {\mathbb N}, i_3 = \max \{ i_0, i_1, i_2 \} $ such that,

\be \int_{\Omega} |\Delta [(u_i-u)\tilde \eta_{\epsilon}]|dx \leq 4 \pi-\epsilon_0+  \epsilon 2 C_1(2 k_1 ||\nabla u||_{L^{\infty}} + 2 + C) \label{(8)}\ee

We choose $ \epsilon >0 $ small enough to have a good estimate of  $ (\ref{(1)}) $.

Indeed, we have:

\begin{displaymath} \left \{ \begin {split} 
      \Delta [(u_i-u) \tilde \eta_{\epsilon}]   & = g_{i,\epsilon}                   \,\, &&\text{in} \!\!&&\Omega \subset {\mathbb R}^2, \\
                 (u_i-u) \tilde \eta_{\epsilon}    & = 0  \,\,                                     && \text{in} \!\!&&\partial \Omega.               
\end {split}\right.
\end{displaymath}

with $ ||g_{i, \epsilon} ||_{L^1(\Omega)} \leq 4 \pi -\epsilon_0/2. $

We can use Theorem 1 of [6] to conclude that there is $ q\geq \tilde q >1 $ such that:

$$ \int_{V_{\epsilon}(x_0)} e^{\tilde q|u_i-u|} dx \leq \int_{\Omega} e^{q|u_i -u| \tilde \eta_{\epsilon}} dx \leq C(\epsilon,\Omega). $$
 
where, $ V_{\epsilon}(x_0) $ is a neighberhooh of $ x_0 $ in $ \bar \Omega $. Here we have used that in a neighborhood of $ x_0 $  by the elliptic estimates, 
$ 1- C \epsilon \leq \tilde \eta_{\epsilon} \leq 1 $. (We can take, $ f(B_{\epsilon^3}(0))$ and we have $ B_{k_2\epsilon^3}(x_0) \subset f(B_{\epsilon^3}(0))\subset B_{k_1\epsilon^3}(x_0) $ for a chart $ (f,B_1(0)) $ around $ x_0 $).

\smallskip

Thus, for each $ x_0 \in \partial \Omega - \{ \bar x_1,\ldots, \bar x_m \} $ there is $ \epsilon_{x_0} >0, q_{x_0} > 1 $ such that:

\be \int_{B(x_0, \epsilon_{x_0})} e^{q_{x_0} u_i} dx \leq C, \,\,\, \forall \,\,\, i. \label{(9)}\ee

By the elliptic estimates (see [13]) $ (u_i)_i $ is uniformly bounded in $ W^{2, q_1}(V_{\epsilon}(x_0)) $ and also, in $ C^1(V_{\epsilon}(x_0)) $. Finaly, we have, for some $ \epsilon > 0 $ small enough,

$$ || u_i||_{C^{1,\theta}[B(x_0, \epsilon)]} \leq c_3 \,\,\, \forall \,\,\, i. $$

We have proved that, there is a finite number of points $ \bar x_1, \ldots, \bar x_m $ such that the squence $ (u_i)_i  $ is locally uniformly bounded (in $ C^{1,\theta}, \theta >0 $) in $ \bar \Omega - \{ \bar x_1, \ldots , \bar x_m \} $.

\underbar {\it Proof of theorem 1.2:} 

We know that:

$$ u_i \in W^{2, k}(\Omega)\cap C^{1, \epsilon}(\bar \Omega). $$ 

We can do integration by parts. The first Pohozaev identity applied around each blow-up point see for example [16] gives :

\be \int_{\partial \Omega_{x_k}} [(\partial_{\nu} u_i) \nabla u_i - \dfrac{1}{2} ||\nabla u_i|^2 \nu] dx =  \int_{\Omega_{x_k}} \nabla V_i e^{u_i} - \int_{\partial \Omega_{x_k}} V_i e^{u_i} \nu, \label{(10)}\ee

Here $ \Omega_{x_k} $ is a neighborhood of $ x_k $ on which we can use the integration by part obtained by a chart around $ x_k $.

\smallskip

We use the boundary condition on $ \Omega $ and the boundedness of $ u_i $ and $ \partial_j u_i $ outside the $ x_k $, to have:

\be \int_{\partial \Omega} (\partial_{\nu} u_i)^2 dx  \leq c_0(b, A, C, \Omega). \label{(11)}\ee

Thus we can use the weak convergence in $ L^2(\partial \Omega) $ to have a subsequence $ \partial_{\nu} u_i $, such that: 

\be \int_{\partial \Omega} \partial_{\nu} u_i \phi dx  \to  \int_{\partial \Omega} \partial_{\nu} u \phi dx, \,\,\, \forall \,\,\, \phi \in L^2(\partial \Omega), \label{(12)}\ee

Thus, $ \alpha_j = 0 $, $ j=1, \ldots, N $ and $ (u_i) $ is uniformly bounded.

\begin{center}
{\bf  ACKNOWLEDGEMENT. }
\end{center}

The work was supported by YY. Li's Grant and Pr Jiguang Bao. The author would like to thank Pr. YY.Li for his support.

\end{document}